
\magnification=\magstep1
\input amstex
\documentstyle{amsppt}
\leftheadtext{K. J. B\"or\"oczky, E. Makai, Jr.}
\rightheadtext{Convex curves enclosing the origin several times}
\topmatter
\title Inverse Blaschke-Santal\'o inequality for convex curves enclosing the
origin several times \endtitle
\author K. J. B\"or\"oczky, E. Makai, Jr.
\vskip.5cm
\centerline{\rm{Alfr\'ed R\'enyi Institute of Mathematics, Hungarian Academy
of Sciences}}
\centerline{H-1364 Budapest, Pf. 127, Hungary}
\vskip.1cm
{\centerline{\rm{http://www.renyi.mta.hu/\~{}carlos,
http://www.renyi.mta.hu/\~{}makai}}}
\vskip.1cm
{\centerline{\rm{E-mail: carlos\@renyi.mta.hu, makai.endre\@renyi.mta.hu}}}
\endauthor
\abstract H. Guggenheimer generalized the planar volume product problem for
locally convex curves $C$ enclosing the origin $k \ge 2$ 
times. He conjectured that the
minimal volume product $V(C)V(C^*)$ for these curves
is attained if the curve consists of the longest
diagonals of a regular $(2k+1)$-gon with centre $0$, 
taken always in the positive orientation.
This conjectured
minimum is of the form $k^2 + O(k)$. We investigate special cases of this
conjecture. We prove it for locally convex
$n$-gons with $2k+1 \le n \le 4k$, if the central
angles at $0$ of all sides are equal to $2k \pi /n$. For $4k+1 \le n$ we prove
that for locally convex $n$-gons enclosing the origin $k \ge 2$ times
the critical (stationary) values of the volume product $V(K)V(K^*)$ 
are attained exactly when up to a non-singular linear map the vertices lie on
the unit circle about $0$, and the central angles of all sides are equal to
$2k \pi /n$. For locally convex
$n$-gons enclosing the origin $k \ge 2$ 
times, and inscribed to the unit circle, with $2k+1 \le n$, 
we prove the conjecture up to a multiplicative factor about $0.43$.
\endabstract
\endtopmatter
\document


{\it 2010 Mathematics Subject Classification.} {\rm{Primary:
52A40. 
Secondary:
52A30, 
52A10


\head
{\bf{1.}} Introduction
\endhead

We begin with some notations, and some well-known facts about the volume
product. Cf., e.g., \cite{L}, \cite{BMMS} and \cite{Mak}.

A {\it{convex body}}\, in ${\Bbb R}^d$ is a compact convex set with nonempty
interior. For $K \subset {\Bbb R}^d$ a convex body with $0 \in 
{\text{int}}\,K$ we write $K^*$ for its polar body, which also is a convex
body, containing $0$ in its interior. One is interested in the infimum and
supremum of $V(K^*)$, if $V(K)$ is given. These 
are of the form $c_{1,d}/V(K)$ and
$c_{2,d}/V(K)$, for 

\newpage

certain constants $c_{i,d}$. (An important special case is
that of $0$-symmetric convex bodies, which however we will not treat.)

We have 
$$
c_{1,d} = \min \{
V(K)V(K^*) \mid K \subset {\Bbb R}^d {\text{ is a convex body with }}0 \in 
{\text{int}}\,K \} .
$$ 
K. Mahler \cite{Mah39} 
conjectured that $c_{1,d} = (d+1)^{d+1}/(d!)^2$, with equality only for a
simplex of barycentre $0$, which is still unproved, although in
many special cases it is known to hold. In particular, for $d=2$ this was
proved by \cite{Mah38}, and the only case of equality is a triangle with
barycentre at $0$, as proved by \cite{Me}. The conjecture is proved, up to a
factor $\left( \pi / (2e) + o(1) \right) ^d$, by \cite{K}.
(For the $0$-symmetric case the analogous minimum is conjectured to be
$4^d/d!$, and is conjectured to be attained e.g. for a parallelotope, or a
cross-polytope, and more generally, it is conjectured to be attained exactly
for those bodies $K$, which, as unit balls of finite dimensional
Banach spaces --- i.e., of Minkowski
spaces --- can be obtained from $[-1,1] \subset {\Bbb{R}}$, 
by taking, in an arbitrary order,
$l^1$-sums and $l^{\infty }$-sums of lower dimensional such Banach spaces.
This conjecture is proved, up to a
factor $\left( \pi / 4 + o(1) \right) ^d$, by \cite{K}. Quite recently the
proof of the three-dimensional case of this conjecture, together with the
conjectured equality cases, i.e., parallelepiped and affine regular
octahedron, was announced in \cite{IS}, and was reassured in \cite{I}.) 

However, $c_{2,d} = \infty $. Therefore one has to consider
the minimax problem, i.e., the supremum of $\min \{ V \left( (K-x)^* \right)
\mid x \in {\text{int}}\,K \} $. (The function $V \left(
(K-x)^* \right) $ is strictly convex for $x \in {\text{int}}\,K$ and tends to
infinity if ${\text{dist}}\,(x, {\text{bd}}\,K) \to 0$, therefore this
function has a unique minimum place, the so called {\it{Santal\'o point $s(K)$
of}}\, $K$.)

This supremum is of the form  $c'_{2,d}/V(K)$, for a certain 
constant $c'_{2,d}$. One has 
$$
\cases
c'_{2,d} = \max \{
V(K) \cdot \min \{ V\left( (K-x)^*) \right) 
\mid x \in {\text{int}}\,K \} \mid K \subset {\Bbb R}^d 
{\text{ is a convex body}} \} \\
= \max \{
V(K)V\left( \left( K-s(K) \right) ^* \right)  \mid K \subset {\Bbb R}^d 
{\text{ is a convex body}} \} .
\endcases
$$
The theorem that $c'_{2,d} = \kappa _d ^2$ (where $\kappa _d$ is the volume of
the unit ball of ${\Bbb R}^d$) was
proved by W. Blaschke and L. Santal\'o, cf. \cite{B} and \cite{S}. The only
case of equality is for the ellipsoid, which was proved by \cite{SR},
\cite{P} and \cite{MP}.

We remark that $V(K)V(K^*)$ is invariant under non-singular linear mappings,
and $V(K)V \left( K - s(K) \right) $ is affine invariant. In particular, if
$K$ admits an affinity with a single fixed point, then this point is the
Santal\'o point of $K$.

These two quantities turn out to be in the cross-road of many disciplines:
they
arose in affine differential geometry (in \cite{B}) and geometry of numbers (in
\cite{M39}), but  also the $0$-symmetric case is very important in finite
dimensional Banach spaces. Further such
disciplines are discrete geometry, geometrical probabilities,
integral geometry in Minkowski spaces, differential equations, and even
theory of functions 

\newpage

of several complex variables. 
 

\head
Inverse Blaschke-Santal\'o inequality for convex curves enclosing the
origin several times
\endhead


H. Guggenheimer \cite{G} posed an interesting generalization of the planar
volume product problem. 

\definition{Definition 1} (H. Guggenheimer, \cite{G}) Let $k \ge 2$ be an
integer.
{\cite{G}} defines the class ${\Cal{C}}_k^2$ of closed curves $C \subset
{\Bbb{R}}^2 \setminus \{ 0 \} $ as follows. $C \in {\Cal{C}}_k^2$ is given 
in polar coordinates as the graph of its radial function  $\varrho ( \cdot )$,
where $\varrho : [0, 2k \pi ] \to (0, \infty )$ (or one can say,
$\varrho ( \cdot )$ is defined on ${\Bbb{R}}$ and is $2k \pi $-periodic).
(Therefore $C$ encircles the origin $k$ times in the positive sense --- in other
words, the winding number, i.e., index of $C$ with respect to $0$ equals $k$.)
Moreover, $C$ is at each of its points locally convex,
that is, has local supporting lines at each of its points, and is
seen from the origin always in the concave side. We topologize  ${\Cal{C}}_k^2$
by the supremum distance of the radial functions, on $[0, 2k \pi ]$ (or on
${\Bbb{R}}$).
\enddefinition


Such a curve also has a 
a {\it{support function $h: [0,2k \pi ] \to (0, \infty )$ (or one can say,
$h( \cdot )$ is defined on ${\Bbb{R}}$ and is
$2k \pi $-periodic)}}. 

One can topologize ${\Cal{C}}_k^2$ also by the 
maximum distance of the support functions. The two definitions are
equivalent: a basic neighbourhood of $C \in {\Cal{C}}_k^2$, with radial function
$\varrho _C(\varphi )$, consists in any of the two cases of those curves
$C' \in {\Cal{C}}_k^2$, whose points $(\varphi , \varrho ')$ satisfy
$(1 - \varepsilon ) \varrho _C (\varphi ) \le \varrho ' \le (1 + \varepsilon )
\varrho _C (\varphi )$, for all $\varphi \in [0, 2k \pi ]$ (or the analogous
inequality for the support functions), where $\varepsilon > 0$ is arbitrary.


The {\it{area $V(C)$ enclosed by $C$}}
can be defined as 
$$
\int _0 ^{2k \pi } \varrho ( \varphi )^2 d \varphi /2
$$ 
(this is also the integral of the index of $C$ on
the whole plane). One can define {\it{polarity on ${\Cal{C}}_k^2$}} 
as usual ({\it{the polar curve is denoted
by $C^*$}}, and has the same properties as $C$), and also the Santal\'o point
as usual, with the usual characteristics. We will use the Santal\'o point only
in the plane.

Observe that 
$$
\cases
{\text{the value of }}
V [ \left( C - x \right) ^*] {\text{ tends to infinity uniformly if }} \\
x {\text{ approaches the closest point of }} C {\text{ on any ray from }} 0. 
\endcases
\tag *
$$
If the ray points to an angle $\varphi \in
[0, 2 \pi )$, then the distance of this closest point is $\min \{ \varrho
(\varphi + 2i \pi ) \mid 0 \le i \le k-1 \} $.
These closest points enclose a convex disc

\newpage

containing $0$ in its interior, we call it the {\it{kernel of $C$}}. 
In particular, the Santal\'o point belongs to the kernel of $C$.

The question is again, as in
the usual case: given $V(C)$, what is the range of values of $V(C^*)$? Again
we have the usual affine invariance property, hence this question is again
equivalent with the question of supremum/infimum, or possibly maximum/minimum 
of the product $V(C) V[\left( C-s(C) \right) ^*]$. (Actually, for this only
dilations with positive ratio would be sufficient.)

\cite{G} also considers the special case of $0$-symmetric curves $C$,
however does not clarify, what does he mean by this. We think that 
the natural way is the
following: {\it{the radial (or support) function is not just $2k \pi
$-periodic, but
actually is $k \pi $-periodic}}. For $k$ even, this means just a curve of
index $k/2$ about $0$, traversed twice, but then the curve $C$ as a set is not
$0$-symmetric in general. Its volume product is $4$ times the volume product
of this curve of index $k/2$ about $0$, hence this case is 
covered by the study of curves of index $k/2$ about $0$, hence is
not to be investigated. 
For $k$ odd, the curve $C$ as a set,
is $0$-symmetric. {\cite{G}} does not seem to recognize these two cases.


For dimension $d \ge 3$ there is an analogous definition, cf. Definition 2
below. For smooth manifolds $M,N$, an immersion $i:M \to N$ 
is a map everywhere 
of rank ${\text{dim}}\,M$. For topological manifolds $M,N$ we say that 
{\it{a map $i:M \to N$ is an immersion}} if $i$ is locally a homeomorphism 
onto its image. That is, each $m \in M$ has a neighbourhood $U$ in $M$ such that
$i$, considered as a map $U \to i(U)$, is a homeomorphism.


\definition{Definition 2} Let $d \ge 3$ and $k \ge 2$ be integers. We write 
${\Cal{C}}^d_k$ for the set of immersed manifolds in ${\Bbb{R}}^d \setminus \{
0 \} $, by an immersion $i : M \to {\Bbb{R}}^d \setminus \{ 0 \} $, where $M$
is a connected compact topological $(d-1)$-manifold, 
for which for any of point $m \in M$ we
have that some open neighbourhood $U \subset M$ of $m \in M$ 
has an image $i(U)$
such that besides $i$ being a homeomorphism $U \to i(U)$, additionally we have
that $i(U)$ is a relatively open subset of the boundary of 
some convex body $K_m$ with $0 \in {\text{int}}\,K_m$. 
\enddefinition


In Definition 2 $i$ is not an embedding: actually, the restriction of $i$
to each linear $2$-subspace of ${\Bbb{R}}^d$ is not an embedding. 
The {\it{kernel of $i(M)$}} 
is the maximal open star domain in ${\Bbb{R}}^d \setminus i(M)$.
The {\it{enclosed volume}} 
is defined as usual, by 
$$
V \left( i(M) \right) = 
V(i,M) = \int _M \| i(m) \| \langle i(m), n \left( i(m) \right) 
dS \left( i(m) \right) /d ,
$$ 

\newpage

where $\| \cdot \|$ is the norm, 
$dS \left( i(m) \right) $ is the surface area measure 
($(d-1)$-Hausdorff measure) element on $i(M)$ at $i(m)$, and
$n \left( i(m) \right) $ is the outer unit normal 
on $i(M)$ at $i(m)$, uniquely defined
$dS \left( i(m) \right) $-almost everywhere. We will use the notation
$V\left( i(M) \right)$, if the immersion $i$ is understood.
Observe that necessariy $M$ is orientable. In fact, $0$ lies always in one of
the open halfspaces bounded by local supporting hyperplanes of 
$i(M)$. The outward normals of these local supporting hyperplanes of 
$i(M)$ will be considered as outward normals of the entire $i(M)$. (These
considerations can be taken over also for $M$ disconnected.)


Below, in the proof of Proposition 3, 
we will see examples of such immersed manifolds. (They will come from mappings 
$S^{d-1} \to S^{d-1}$ of index $k$.)
The significance of
connectedness of $M$ will be explained later, in Remark 8.


For the planar case,
{\cite{G}} tries to apply some local arguments for the lower
estimate. However, in absence of compactness, local arguments are definitely 
insufficient for this. And in fact, for each $k \ge 2$, 
the equivalence classes of the
curves in ${\Cal{C}}^2_k$ with respect to nonsingular linear maps 
do not form a compact set in their natural topology,
contrary to what is asserted in {\cite{G}}. More exactly, we have Corollary 4
below.


\proclaim{Proposition 3}
Let $d, k \ge 2$ be integers. Then the continuous
affine invariant functional $V(C) V[ \left( C-s(C) \right) ^*]$ is unbounded
above. In other words, there is no Blaschke-Santal\'o theorem for
${\Cal{C}}_k^d$.
\endproclaim


\demo{Proof}
First let $d=2$.
We write, as usual, $e_1=(1,0)$ and $e_2=(0,1)$. Let $\varepsilon > 0$ be an
arbitrarily small number. We define $C$ as follows. Let $E_1$ and $E_2$ be
ellipses with equations $x^2 + (y/\varepsilon )^2 = 1$ and 
$x^2 + (y \varepsilon )^2 = 1$. We define $C$ as follows. First we traverse
${\text{bd}}\,E_1$ once, from $e_1$ to $e_1$, in the positive sense, and then in
continuation we traverse ${\text{bd}}\,E_2$, $k-1$ times, from 
$e_1$ to $e_1$, in the positive sense. Let $x \in {\text{int}}\,(E_1 \cap E_2)
= {\text{int}}\,E_1$ (this set is the kernel of $C$). Then
$$ 
V(C) \ge V(E_2) = \pi / \varepsilon , {\text{ and }} 
V \left( (C-x)^* \right) \ge 
V \left( (E_1-x)^* \right) \ge V(E_1^*) = \pi / \varepsilon ,
$$
hence 
$$
\cases
V(C) V[ \left( C - s(C) \right) ^* ] = \\
V(C) \inf \{ V \left( (C-x)^* \right) \mid x \in 
{\text{\rm{kernel of }}} i(M) \} \ge \pi ^2/ \varepsilon ^2 .
\endcases
$$
 
For $d=3$ we rotate the $2$-dimensional example about the $x_1$-axis, and
in the same way 
we obtain the $(d+1)$-dimensional example from the $d$-dimensional
example. Then we obtain, writing $\kappa _d$ for the volume of the unit ball in
${\Bbb{R}}^d$, that
$$
\cases
V(C) V[ \left( C - s(C) \right) ^* ] = \\
V(C) \inf \{ V \left( (C-x)^* \right) \mid x \in 
{\text{\rm{kernel of }}} i(M) \} \ge \kappa _d ^2 / \varepsilon ^d .
\endcases
$$


\newpage

\proclaim{Corollary 4}
Let $d=2$, and let $k \ge 2$ be an integer. Then the 
equivalence classes of the immersed manifolds $C \in {\Cal{C}}_k^2$ with
respect to non-singular linear maps
do not form a compact set in the quotient topology.
\endproclaim


\demo{Proof}
By \thetag{*}
the product $V(C) V \left( C - s(C) \right) =
V(C) \inf \{ V \left( (C-x)^* \right) \mid x \in {\text{kernel of }} 
i(M) \} $ is continuous, is invariant with respect to non-singular
linear maps, and is also unbounded above, by Proposition 3. Therefore
the equivalence classes of 
the immersed manifolds from ${\Cal{C}}^2_k$ with respect to non-singular
linear maps do not form a compact set in their
quotient topology. 
$\blacksquare $
\enddemo


Thus the remaining question is whether there is an inverse Blaschke-Santal\'o
inequality here. Below we will
describe the (conjecturable) statement of \cite{G}, for the planar case. 
We begin with a notation.


\definition{Notation 5}
Let $k \ge 2$ and $n \ge 2k+1 \ge 5$. Then $C_{n,k}$ 
denotes the following
closed polygonal line. We consider a regular $n$-gon of centre $0$, inscribed
in the unit circle $S^1$ of centre $0$,
and we pass successively on its
$k$-th smallest diagonals, always in the positive sense, until the sum of the
central angles of the sides attains $2k \pi $. We write $(k,n)$ for the
greatest common divisor of $k,n$. Then
$C_{n, k}$ 
passes only on each $(k,n)$-th vertex of the regular $n$-gon,
but passes through each of them $(k,n)$ times.
\enddefinition


\definition{Conjecture 6} (H. Guggenheimer {\cite{G}}, stated there as theorems)
1) For the general (i.e., 
not $k \pi $-periodic) case {\cite{G}} considers $C_{2k+1,k}$.
This has index $k$
with respect to $0$, and is conjectured to give the minimal volume product 
$V(C)V(C^*)$ for all $C \in {\Cal{C}}_k^2$.

2) For the $k \pi $-periodic case {\cite{G}} considers $C_{2k+2,k}$.
This has index $k$ with respect to $0$, and is conjectured 
to give the minimal volume product $V(C)V(C^*)$ for the $k \pi $-periodic
case. (We have to remark that for $k$ even we have here a doubly traversed 
curve $C_0 \in {\Cal{C}}_{k/2}$, and then $V(C)V(C^*) = 4 V(C_0)V(C_0^*)$.
Therefore part 1) of this conjecture, for index $k/2$, implies part 2) of this
conjecture for index $k$. However, for $k$ odd, the conjectured $C_{2k+2,k}$ 
passes through all vertices of the regular $(2k+2)$-gon, through each of them
just once, and then
part 2) of this conjecture does not follow from its part 1). Moreover, for $k$
odd, the curve $C_{2k+2,k}$ is $0$-symmetric.)
\enddefinition


\definition{Remark 7}
For $C = C_{n,k}$, 
with $n \ge 2k+1$, we have, writing $\vartheta _i$ for the
central angle of the $i$'th side, that, $\sum _{i=1} ^n \vartheta _i =
2k \pi $, all $\vartheta _i$'s equal $2k \pi /n$, 
and $V(C)V(C^*) = [ \sum _{i=1}^n ( \sin \vartheta _i ) /2 ]
\cdot [ \sum _{i=1} ^n \tan (\vartheta _i/2) ] = 
n^2 \sin ^2 (k \pi /n) = \big[ [
\left( \sin (k \pi /n) \right) / ( k \pi /n) ] \cdot k \pi \big] ^2$, which
strictly increases with $n$. So for a minimum we must have the minimal
possible $n$, i.e., $n = 2k+1$ for the general case, and $n=2k+2$ for 
the $k \pi $-periodic case with $k$ odd. 
This is a small support for our conjecture.
The respective values of $V(C)V(C^*)$ in these two cases
are $(2k+1)^2 \sin ^2 \left( k \pi / (2k+1) \right) $, and 
$(2k+2)^2 \sin ^2 \left( k \pi / (2k+2) \right) $, both of the form $4k^2
+O(k)$. 
\enddefinition


\newpage

\definition{Remark 8}
Suppose that 
$M$ has several connected components $M_1, \ldots ,M_l$ (fi\-ni\-te\-ly many by
compactness of $M$), with respective indices, i.e., winding numbers, 
$k_1, \ldots , k_l$, satisfying
$0 < k_j$ and $\sum _{j=1}^l k_j = k$. 
Then the inverse
Blaschke-Santal\'o inequality can be asked also for $M$. However, this
question
can be reduced to the connected case, with smaller indices $k$. In fact, we
have $V(iM) = \sum _{j=1}^l V(iM_j)$ and 
similarly $V \left( (iM)^* \right) = \sum _{j=1}^l V \left( (iM_j)^* \right)
$, hence by the arithmetic-geometrical mean inequality
$$
\cases
V(iM)V \left( (iM)^* \right) \ge 
\left( \prod _{j=1}^l V(iM_j) \right) ^{1/l} \cdot \\
\left( \prod _{j=1}^l V \left( (iM_j)^* \right) \right) ^{1/l} = 
\left( \prod _{j=1}^l V(iM_j) V \left( (iM_j)^* \right) \right) ^{1/l} .
\endcases
$$
Hence if we have some non-trivial lower
estimates for the indices $k_1, \ldots , k_l$, then this implies some
nontrivial lower estimate for $V(iM)V \left( (iM)^* \right) $. For the planar
case, assuming that Conjecture 6 were
valid for indices smaller than $k$, we would have
Conjecture 6 for disconnected $M$ with index $k$.
\enddefinition


We turn to the case $d=2$.
By an approximation argument, it is
sufficient to prove this conjecture
for $n$-gons, where $n \ge 3$ is an arbitrary 
integer.
We use polar coordinates, i.e., the vertices will be given as $(\varphi
_i, \varrho _i)$, where $\varphi _i \in [0, 2k \pi ]$ --- or $\varphi _i \in
{\Bbb{R}}$ and $\varrho (\varphi )$ is $2k \pi $-periodic --- 
and $\varrho _i \in (0, \infty )$.
Since the central angles of the polygon are less than $\pi $, and
their sum is $2k \pi $, therefore necessarily we have $n \ge 2k+1$. 

Some numerical experimentation suggests for the case $k=2$ that in part 1) of
the Conjecture we have actually a local minimum among pentagons, 
and for the case $k=3$ that in part 2) of the Conjecture we have actually a 
local minimum among octagons.


\proclaim{Theorem 9}
Let $k \ge 2$ and $n \ge 2k+1$, and let $C \in {\Cal{C}}_k^2$ be 
a closed $n$-gon
(degeneration to a polygon with less than $n$ vertices is excluded). Then
we have $(1) \Longrightarrow (2)$, where

(1) $C$ is affinely equivalent to $C_{n,k}$,

(2) the volume product $V(C)V[ \left( C - s(C) \right) ^*]$ is a
critical value. 

For $n \ge 4k+1$ we have also $(2) \Longrightarrow (1)$.
\endproclaim


\newpage

\demo{Proof}
We begin with the proof of $(1) \Longrightarrow (2)$.
For $C = C_{n,k}$ with $n \ge 2k+1$
an easy calculation shows that the partial derivatives
of $V(C)V(C^*)$ with respect to the angular and radial coordinates of the
vertices are $0$. Observe that $C_{n,k}$ 
has an $n/(n,k)$-fold rotational
symmetry (and $n/(n,k) \ge n/k \ge (2k+1)/k > 2$), 
hence its Santal\'o point is $0$.
Then applying the statement about the stability of the
Santal\'o point, in \cite{BMMR}, Lemma 11 and \cite{BM}, Theorem E
(valid also for ${\Cal{C}}^2_k$),
we obtain statement (2) for $V(C)V[ \left( C - s(C) \right) ^*]$ 
rather than for $V(C)V(C^*)$.

We turn to the proof of $(2) \Longrightarrow (1)$ for $n \ge 4k+1$.
Observe that the average central angle of the sides
is $2k \pi /n \le 2k \pi /(4k+1) < \pi /2$. Hence the average sum of the
central angles of two adjacent sides is at most $4k \pi / (4k + 1) < \pi
$. Therefore the sum of the
central angles of some two adjacent sides is less than $\pi $.
Using this, our proof for the usual case, i.e., for $k=1$, 
cf. \cite{BM}, Theorems A and F, gives that if $C$ gives
a critical value of $V(C)V[ \left( C - s(C) \right) ^*]$, then some affine
image of $C$ is inscribed to the unit circle about $0$, has a positive
orientation, and has equal sides.
However, this polygonal line must close after $n$ steps, and just after a
total angle of rotation $2k \pi $, hence it is $C_{n,k}$. $\blacksquare $
\enddemo

%
%

We can support our Conjecture 6 by investigating two special cases of it.
Considering $C_{n,k}$ as inscribed to $S^1$, 
we can preserve the angular coordinates of the vertices 
of the conjectured $C_{n,k}$ while changing their radial coordinates,
or we can preserve the radial coordinates of the vertices
of the conjectured $C_{n,k}$
while changing their angular coordinates (and also their number).
For the inverse Blaschke-Santal\'o inequality we give the exact lower bound 
of the volume product $V(C)V(C^*)$ in
the first case, for $2k + 1 \le n \le 4n$, 
and some positive bound in the second case, for $2k + 1 \le n$. 


\proclaim{Proposition 10}
Let \,$2k+1 \le n \le 4k$.
Let $C \in {\Cal{C}}_k^2$ be a closed $n$-gonal
line with vertices having angular coordinates $2 \pi ik/n \in [0, 2k \pi ]$, 
for $0 \le i \le n$ (the $0$'th and $n$-th vertices coincide). Then
$V(C)V(C^*) \ge n^2 \sin ^2 (k \pi /n)$, with equality only if either $C$ is a
copy of $C_{n,k}$, magnified from the origin, or we have $n=4k$ and $C$ is a
$k$ times traversed rhomb of centre $0$.
\endproclaim


\demo{Proof}
We write $\varrho ( \cdot ): [0, 2k \pi ] \to (0, \infty ) $ 
for the radial function of $C$. Further, we write $\varphi _i := 2 \pi ik/n$,
and $\varrho _i :=  \varrho (\varphi _i)$. 
By $\varphi _i = i 2k \pi /n$ the central angle of the side
$[( \varphi _i, \varrho _i ), ( \varphi _{i+1}, \varrho _{i+1}) ]$ is
$ \vartheta _i =  \vartheta := 2k \pi /n$.
Then

\newpage

$$
V(C) = \sum _{i=1}^n \varrho _i \varrho _{i+1} (\sin \vartheta _i) /2.
$$

Now we are going to determine $V(C^*)$. The angular domains, with vertices at
$0$ and with boundary rays passing through the vertices of $C$, decompose
also $C^*$ into $n$ domains. Each of these domains $Q_i$ is a convex
quadrangle, with one vertex at $0$, and two sides beginning at $0$,
lying on the two boundary rays of the angular domain 
$( \varphi _i, \varrho _i ) 0 ( \varphi _{i+1}, \varrho _{i+1} ) $. 
These two sides have lengths $1/\varrho
_i$ and $1/\varrho _{i+1}$, and the other endpoints
of these two sides have right angles in $Q_i$. The diagonal of $Q_i$ from $0$
decomposes $Q_i$ into two right triangles, and $V(C^*)$ is the sum of the 
areas of these two triangles. We may suppose that the angle bisector of
the angle of $Q_i$ at $0$ is the positive $x$-axis. Then 
we obtain by a elementary calculation that 
$$
\cases
4V(Q_i) = [1/ \varrho _i ^2 + 1/ \varrho _{i+1}^2] \cdot
[-\cos ^2 (\vartheta /2) \cot (\vartheta /2) + 
\sin ^2 (\vartheta /2) \tan (\vartheta /2)] \\
+[2 \cdot (1/\varrho _i) \cdot  (1/ \varrho _{i+1})] \cdot 
[\cos ^2 (\vartheta /2) \cot
(\vartheta /2) + \sin ^2 (\vartheta /2) \tan (\vartheta /2) + \sin \vartheta ].
\endcases
$$
We write $F(\vartheta )$ and $G(\vartheta )$ for the coefficients of
$1/ \varrho _i ^2 + 1/ \varrho _{i+1}^2$ and $2 \cdot (1/\varrho _i) \cdot 
(1/ \varrho _{i+1})$ in this formula. 
Then $F(\vartheta ) \ge 0$: more exactly, for $n=4k$ we have $F(\vartheta ) =
0$, 
and for $2k +1 \le n \le 4k-1$ we have $F(\vartheta ) > 0$. In fact, this
last equality and inequality follow from
$$
\cases
[\cos ^2 (\vartheta /2) \cot (\vartheta /2)] / 
[\sin ^2 (\vartheta /2) \tan (\vartheta /2)] = \cot ^4 (\vartheta /2) = \\
\cot ^4 \left( (2k \pi /n) /2 \right) \le 
\cot ^4 \left( \left( 2k \pi /(4k) \right) /2 \right) = \cot ^4 (\pi /4) = 1,
\endcases
$$
and in the inequality here we have equality for $n=4k$, and strict inequality
for $2k+1 \le n \le 4k-1$. 

By all these calculations we have
$$
\cases
V(C)V(C^*) =
[\sum _{i=1}^n \varrho _i \varrho _{i+1} (\sin \vartheta ) /2] \cdot \\
\sum _{i=1}^n \left[
(1/ \varrho _i ^2 + 1/ \varrho _{i+1}^2) \cdot F(\vartheta )/4 +
\left( 2/(\varrho _i \varrho _{i+1}) \right) \cdot G(\vartheta )/4 \right].
\endcases
\tag A
$$
Here the indices are considered cyclically, 
for both sums here, i.e., for $S_1 := V(C)$ and $S_2 := V(C^*)$. 
Applying the arithmetic-geometric mean inequality both for $S_1$ and $S_2$, 
from \thetag{A} we obtain
$$ 
\cases
V(C)V(C^*) \ge n \cdot [\prod _{i=1}^n (\varrho _i \varrho _{i+1}) ]^{1/n} 
[(\sin \vartheta ) /2] \cdot \\
n \cdot \prod _{i=1}^n \left[
(1/ \varrho _i ^2 + 1/ \varrho _{i+1}^2) \cdot F(\vartheta )/4 + 
\left( 2/(\varrho _i \varrho _{i+1}) \right) 
\cdot G(\vartheta )/4 \right] ^{1/n} 
= \\
n^2 \cdot [(\sin \vartheta ) /2] \cdot \prod _{i=1}^n \left[
(\varrho _{i+1} / \varrho _i + \varrho _i / \varrho _{i+1}) 
\cdot F(\vartheta )/4 + 2 \cdot G(\vartheta )/4 \right] ^{1/n} .
\endcases
\tag B
$$

\newpage

In \thetag{B}, in the arithmetic-geometric mean inequality for $S_1$ 
we have equality, if and only if all its summands are equal, i.e., if 
for each $i$ we have $\varrho _i \varrho _{i+1} =
\varrho _{i+1} \varrho _{i+2}$, i.e., 
$$
\varrho _i = \varrho _{i+2}
\tag C
$$ 
(since $ 
\varrho _{i+1} > 0$). That is, for $n$ odd all $\varrho _i$'s are equal, in
which case the statement of the theorem is proved, while
for $n$ even the $\varrho _i$'s with $i$ of given parity are equal  --- i.e.,
the $\varrho _i$'s assume alternately two values. (The analogous consideration
for $S_2$ will not be needed.)

In \thetag{B}, under the last product sign, 
again  by the arithmetic-geometric mean inequality, we have for each $i$
that
$$
(\varrho _{i+1} / \varrho _i + \varrho _i / \varrho _{i+1}) 
\cdot F(\vartheta )/4 \ge 2 \cdot F(\vartheta )/4 ,
\tag D
$$
with equality for any $i$ if only if 
$$
\cases
{\text{either }} F(\vartheta ) = 0, {\text{ i.e., }} n = 4k, 
{\text{ or }} F(\vartheta ) > 0, {\text{ i.e.,}} \\
2k+1 \le n \le 4k-1, 
{\text{ and }} \varrho _i = \varrho _{i+1} {\text{ for each }} i .
\endcases
\tag E
$$ 

From \thetag{B} and \thetag{D} we obtain
$$  
\cases
V(C)V(C^*) \ge 
n^2 \cdot [(\sin \vartheta ) /2] \cdot \prod _{i=1}^n (2 \cdot F(\vartheta )/4
+ 2 \cdot G(\vartheta )/4 ) ^{1/n} 
= \\
n^2 \cdot [(\sin \vartheta ) /2] \cdot 
\left( F(\vartheta ) + G(\vartheta ) \right) /2 
= \\
n^2 \cdot [(\sin \vartheta ) /2] \cdot \sin
(\vartheta /2) \left( \cos (\vartheta /2) + \sin (\vartheta /2) 
\tan (\vartheta /2) \right) =  \\
n^2 \sin ^2 (\vartheta /2) = n^2 \sin ^2 (k \pi /n) ,
\endcases
\tag F
$$
with equality only if both \thetag{C} and \thetag{E} hold.
In other words, either $n = 4k$, and 
$\varrho _i$'s assume alternately two values, or $2k+1 \le n \le 4k-1$ and
all $\varrho _i$'s are equal.
In other words, we have $n=4k$ and $C$ is a $k$ times traversed rhomb of
centre $0$, thus is a non-singular linear image of $C_{n,k} = C_{4k,k}$ (in
particular, $V(C)V(C^*) = V(C_{4k,k})V(C_{4k,k}^*)$), or
$2k+1 \le n \le 4k-1$ and all $\varrho _i$'s are equal, i.e., we have that $C$
is an inflation from $0$ of $C_{n,k}$.
$\blacksquare $
\enddemo


The following Proposition 12 proves Conjecture 6 in a special case, up to a
constant factor about $0.43$. Before it we need a lemma.


\proclaim{Lemma 11}
The functions $1/(1 - \cos t)$ and $t/ \sin t$ are strictly
convex for $t \in (0, \pi )$.
\endproclaim


\newpage

\demo{Proof}
We begin with $1/(1-\cos t)$. Its derivative is 
$$
- \sin t / (1 - \cos t)^2 = - \cos (t/2) /[2 \sin ^3 (t/2)], 
$$
which is strictly
increasing since $\cos (t/2)$ is strictly decreasing and $\sin (t/2)$ is
strictly increasing for $t \in (0, \pi )$.

Next we deal with $t/ \sin t$. Its second derivative is
$$
(t+t \cos ^2 t - 2 \cos t \sin t)/\sin ^3 t,
$$
and we have to show that here the numerator is positive. Equivalently,
writing $s := 2t$,
$$
t > 2 \cos t \sin t / (1 + \cos ^2 t), {\text{ i.e., }} s \ge (4 \sin
s)/(3 + \cos s),
$$
for $s = 2t \in (0, 2 \pi )$. 
However, the left/right hand side of the last inequality is 
positive/non-positive for 
$s \in [ \pi ,  2 \pi )$, therefore the strict inequality holds here.
Hence we suppose $s \in (0, \pi )$. Both sides of the
last inequality are $0$ for $s=0$, so it suffices to prove the respective
inequality for the derivatives of the left and right hand side expressions.
I.e., we have to prove 
$$
1 > (12 \cos s + 4) / ( 3 + \cos s)^2  .
$$
Rearranging, this becomes
$$
(5 - \cos t) (1 - \cos t) > 0, 
$$
which is valid for $s \in (0, \pi )$.
$\blacksquare $
\enddemo


\proclaim{Proposition 12}
Let $k \ge 2$ and $2k+1 \le n$ be integers.
Let $C \in {\Cal{C}}_k^2$ be a closed polygonal
line inscribed in the unit circle about $0$. Then
$V(C)V(C^*) \ge 4k^2/ \min \{ (4 / \pi ^2) (t / \sin t) + 2 / (1 - \cos t) \mid
t \in (0, \pi ) \} = k^2 \cdot 1.7366\ldots $\,. 
\endproclaim


\demo{Proof}
Let $C$ have $n$ vertices.
Let the angles with vertex $0$, spanned by the sides of the closed $n$-gonal
line $C$,
be $\vartheta _i \in (0, \pi )$, for $1 \le i \le n$, where 
$$
\sum _{i=1}^n \vartheta _i = 2k \pi .
\tag A
$$ 
Then (for the first formula cf. the proof of Proposition 10), 
$$
V(C) = \sum _{i=1}^n (\sin \vartheta _i) / 2, 
{\text{ and }} V(C^*) = \sum _{i=1}^n \tan (\vartheta _i/2).
\tag B
$$

\newpage

We choose some constant $c \in (0, \pi )$ (later we will optimize its value).
Then
$$
\cases
V(C)C(C^*) \ge \sum \{ (\sin \vartheta _i) / 2 \mid \vartheta _i \in (0,c] \}
\cdot \sum \{ \tan (\vartheta _i / 2) \mid \vartheta _i \in (0,c] \} + \\
\sum \{ (\sin \vartheta _i) / 2 \mid \vartheta _i \in (c, \pi ) \}
\cdot \sum \{ \tan (\vartheta _i / 2) \mid \vartheta _i \in (c, \pi ) \} .
\endcases
\tag C
$$
We denote 
$$
k_1 := \sum \{ \vartheta _i \mid \vartheta _i \in (0, c] \big/ \pi 
{\text{ and }}
k_2 := \sum \{ \vartheta _i \mid \vartheta _i \in (c, \pi ) \big/ \pi ,
\tag D
$$
where (cf. \thetag{A})
$$
k_1,k_2 \in [0,2k], {\text{ and }} k_1 + k_2 = 2k.
\tag E
$$

For the first summand in \thetag{C} we use the estimates 
$$
(\sin \vartheta _i)/2 \ge (\vartheta _i /2) \cdot (\sin c)/c {\text{ and }}
\tan (\vartheta _i /2) \ge \vartheta _i /2.
\tag F
$$
For the second summand in \thetag{C} we use the estimates 
$$
\cases
(\sin \vartheta _i)/2 = \left( \sin (\pi - \vartheta _i) \right) /2 \ge  
[(\pi - \vartheta _i)/2] \cdot 
[\left( \sin ( \pi - c) \right) /( \pi - c)] \\
{\text{and }}
\tan (\vartheta _i /2) = 1 / \tan \left( (\pi - \vartheta _i ) /2 \right) 
\ge \\ 
[1 / \left( (\pi - \vartheta _i ) /2 \right) ] \cdot [ \left( \pi - c )/2
\right) / \tan \left( (\pi - c )/2 \right) ].
\endcases
\tag G
$$

Thus for the first summand in \thetag{C} we obtain by \thetag{F}
$$
\cases
\sum \{ \sin \vartheta _i) / 2 \mid \vartheta _i \in (0,c] \}
\cdot \sum \{ \tan (\vartheta _i / 2) \mid \vartheta _i \in (0,c] \} \ge \\
\sum \{ \vartheta _i / 2 \mid \vartheta _i \in (0,c] \} \left( (\sin c)/c 
\right) \cdot
\sum \{ \vartheta _i / 2 \mid \vartheta _i \in (0,c] \} = \\
(k_1 \pi /2) \left(
(\sin c) / c \right) \cdot (k_1 \pi /2) = k_1^2 \cdot \pi ^2 /4 \cdot
(\sin c) / c .
\endcases
\tag H
$$
Analogously, for the second summand in \thetag{C} we obtain by \thetag{G}, and
by the arith\-me\-tic-harmonic mean inequality that
$$
\cases
\sum \{ \sin \vartheta _i) / 2 \mid \vartheta _i \in (c, \pi ) \}
\cdot \sum \{ \tan (\vartheta _i / 2) \mid \vartheta _i \in (c, \pi) \} 
\ge \\
| \{ \vartheta _i \mid \vartheta \in (c, \pi ) \} |^2 \cdot [\left( \sin ( \pi
- c) \right) / (\pi - c)] \cdot [\left( ( \pi - c)/2 \right) / \tan 
\left( ( \pi - c)/2 \right) ] = \\
| \{ \vartheta _i \mid \vartheta \in (c, \pi ) \} |^2 \cdot (1 - \cos c)/2.
\endcases
\tag I
$$
Since each $\vartheta _i\in (c, \pi )$ is smaller than $\pi $, and their sum is
$k_2 \pi $, therefore for their number we obtain

\newpage

$$
| \{ \vartheta _i \mid \vartheta \in (c, \pi ) \} | \ge k_2,
\tag J
$$
hence \thetag{I} gives
$$
\cases
\sum \{ \sin \vartheta _i) / 2 \mid \vartheta _i \in (c, \pi ) \}
\cdot \sum \{ \tan (\vartheta _i / 2) \mid \vartheta _i \in (c, \pi) \} 
\ge \\
k_2^2 \cdot (1 - \cos c)/2.
\endcases
\tag K
$$
By \thetag{C}, \thetag{H} and \thetag{K} we have
$$
\cases
V(C)C(C^*) \ge \sum \{ \sin \vartheta _i) / 2 \mid \vartheta _i \in (0,c] \}
\cdot \sum \{ \tan (\vartheta _i / 2) \mid \vartheta _i \in (0,c] \} + \\
\sum \{ \sin \vartheta _i) / 2 \mid \vartheta _i \in (c, \pi ) \}
\cdot \sum \{ \tan (\vartheta _i / 2) \mid \vartheta _i \in (c, \pi ) \}
\ge \\
k_1^2 \cdot ( \pi ^2/4) \cdot \left( (\sin c)/c \right) +
k_2^2 (1 - \cos c)/2.
\endcases
\tag L
$$
(Observe that this holds also for $k_1 = 0$ and for $k_2 = 0$.)

Using \thetag{E}, we have $k_1, k_2 \in [0, 2k]$, and 
we substitute $2k - k_1$ for $k_2$ in the last expression in
\thetag{L}. Thus we obtain a
quadratic polynomial $p(k_1)$ of $k_1$. 
Next we minimize the value of $p(k_1)$
for all $k_1 \in [0, 2k]$. Then clearly \thetag{L} will remain
valid if we replace $p(k_1)$ in the last expression of \thetag{L} 
by $\min \{ p(k_1) \mid k_1 \in [0, 2k] \} $.
This minimum is attained for 
$$
k_1 = [2k (1
- \cos c)/2] / [ ( \pi ^2 /4) \left( (\sin c)/c \right) + (1 - \cos c)/2 ]
\in (0, 2k),
\tag M
$$ 
and its value is 
$$
4k^2 \cdot \frac{(\pi ^2/4) \cdot \left( (\sin c) /c \right) \cdot
(1 - \cos c)/2}
{(\pi ^2/4) \cdot \left( (\sin c) /c \right) + (1 - \cos c)/2}.
\tag N
$$
This value still depends on the arbitrarily chosen value $c \in (0, \pi )$. We
have to choose $c$ so that \thetag{N} becomes maximum. Equivalently, we
want to maximize the coefficient of $4k^2$ in \thetag{N}. It will be more
convenient to consider the reciprocal of this coefficient, and then we will
have to look for the minimum of this reciprocal. This reciprocal is
$$
\frac{4c}{ \pi ^2 \sin c } + \frac{2}{1 - \cos c}.
\tag O
$$

By Lemma 11 \thetag{O} is a strictly
convex function of $c \in (0, \pi )$. Moreover, it has
limits at $0$ and $ \pi $ equal to $\infty $. Hence it has a unique minimum
place $c_0$, which is the unique root of its derivative.
That is, 

\newpage

$$
\frac{4}{ \pi ^2} \cdot \frac{ \sin c_0 - c_0 \cos c_0}{\sin ^2 c_0} -
\frac{2\sin c_0}{(1 - \cos c_0)^2 } = 0.
\tag P
$$
Solving this numerically, we find $c_0 \approx 115.5 ^{\circ }$ and 
and the maximum value of \thetag{N} is $k^2 \cdot 1.7366 \ldots $, as asserted
in the Proposition.
$\blacksquare $
\enddemo


\definition{Remark 12}
Analogously to the planar case, possibly for each
fixed $d$, for $C \in {\Cal{C}}^d_k$,
there would hold a lower bound $V(C)V(C^*) \ge {\text{const}}_d \cdot k^2$? 
\enddefinition


\definition{Remark 13}
For $d \ge 3$ and $k \ge 2$
possibly some polyhedral surfaces would minimize $V(C)V(C^*)$ in 
${\Cal{C}}^d_k$. For $\Sigma $ a simplex of barycentre $0$, suitable maps 
${\text{bd}}\,\Sigma \to {\text{bd}}\,\Sigma $ of index $k$ (analogous to 
maps $S^{d-1} \to S^{d-1}$ of index $k$) give $V(C)V(C^*) =
k^2 (d+1)^{d+1}/(d!)^2$. In particular, for $d=3$ they give $(64/9)k^2$.
However, a right prism of height $2$ over the planar conjectured extremal
curve, realized as a direct product with $[-1,1]$, 
with polar the respective bipyramid of height $2$, gives for $k \ge 3$ 
better, while for $k=2$  the doubly traversed simplex gives better.
In fact, by Remark 7, this inequality is 
$$
(4/3)(2k+1)^2 \sin ^2 \left( k \pi / (2k+1) \right) < (64/9)k^2.
$$
This can be rewritten as 
$$
\sin ^2 \left( k \pi / (2k+1) \right) < 16k^2 / [3(2k+1)^2].
$$
Here the right hand side is greater than $1$ for all $k \ge 4$, so in this
case this inequality holds. For $k=3$ a direct calculation shows the same.
However, for $k=2$ the converse inequality holds. 
For higher dimensions, one could take products of lower dimensional
examples, with the indices being multiplied (with the $l^{\infty }$- or the
$l^1$-norm), iteratedly. However, we do not have a
reasonable conjecture for $C \in {\Cal{C}}^d_k$, for any given $d \ge 3$ and
$k \ge 2$, except possibly for $d=3$ and $k=2$ the doubly traversed simplex?
\enddefinition


\Refs

\widestnumber\key{WWW}


[B]
Blaschke, W.,
\"Uber affine Geometrie VII: Neue Extremeigenschaften von Ellipse und
Ellipsoid,
{\it{Ber. \"uber die Verhandl. der K\"onigl. S\"achs. Gesellschaft der
Wiss. zu Leipzig, Math.-Phys. Klasse}}
{\bf{69}} (1917), 306-318, Jahresberichte Fortschr. Math. {\bf{46.}}1112

\newpage

[BM]
B\"or\"oczky, K. J., Makai, E. Jr., {\it{Remarks on planar Blaschke-Santal\'o
inequality}}, arXiv:
\newline
1411.3842

[BMMR]
B\"or\"oczky, K. J., Makai, E. Jr., M. Meyer, S. Reisner,
Volume product in the plane -- lower estimates with stability,
{\it{Stud. Sci. Math. Hungar.}} {\bf{50}} (2013), 159-198,
MR {\bf{3187810}}

[I]
Iriyeh, H., {\it{On Mahler's conjecture -- a symplectic aspect --}},
www.f.waseda.jp/martin/conf/
\newline
dgde/iriyeh.pdf, Dec. 13, 2017

[IS]
Iriyeh, H., Shibata, M., {\it{Symmetric Mahler conjecture for the volume 
product in the three diemnsional case}}, arXiv:1706.01749v2

[G] 
Guggenheimer, H.,
Hill equations with coexisting periodic solutions,
{\it{J. Diff. Equ.}} {\bf{5}} (1969), 159-166, 
MR {\bf{39\#}}550 

[K] Kuperberg. G.,
From the Mahler conjecture to Gauss linking integrals,
{\it{Geom. Funct. Anal.}} {\bf{18}} (2008), 870-892, arXiv:math/0610.5904,
MR {\bf{2009i:}}52005

[L] Lutwak, E., 
Selected affine isoperimetric inequalities,
In: {\it{Handbook of Convex Geometry}} (Eds. P. M. Gruber, J. M. Wills), 
North Holland, Amsterdam etc. 1993, 151-176, MR {\bf{94h:}}{\rm{52014}}

[Mah38] Mahler, K.,
Ein Minimalproblem f\"ur konvexe Polygone,
{\it{Mathematika B (Zutphen)}} {\bf{7}} 
(1938), 118-127, Zbl. {\bf{20.}}{\rm{50}}

[Mah39] Mahler, K.,
Ein \"Ubertragungsprinzip f\"ur konvexe K\"orper,
{\it{\v Casopis P\v est. Mat. Fys.}} {\bf{68}} 
(1939), 93-102, MR {\bf{1,}}{\rm{202}} 

[Mak]
Makai, Jr., E.,
The recent status of the volume product problem, {\it{\'Etudes Op\'eratorielles,
Banach Center Publications}}, {\bf{112}}, Inst. of Math., Polish Acad. Sci.,
Warszawa, 2017, 273-280, MR {\bf{3754082}}

[Me]
Meyer, M.,
Convex bodies with minimal volume product in ${\Bbb R}^2$,
{\it{Monatsh. Math.}} {\bf{112}} (1991), 297-301, MR {\bf{92k:}}{\rm{52015}}

[MP] Meyer, M., Pajor, A.,
On the Blaschke-Santal\'o inequality,
{\it{Arch. Math. (Basel)}}\,\, {\bf{55}} (1990), 
82-93, MR {\bf{92b:}}{\rm{52013}}

[P]
Petty, C. M.,
Affine isoperimetric problems,
{\it{Ann. N. Y. Acad. Sci.}} {\bf{440}} 
(1985), 113-127, MR {\bf{87a:}}{\rm{52014}}

[SR]
Saint Raymond, J.,
Sur le volume des corps convexes sym\'etriques,
{\it{S\'em. d'Initiation \`a l'Analyse 20$^{\text{e}}$ Ann\'ee}}, 
1980-1981, Exp. {\bf{11}} (Univ. Paris VI, Paris, 1981),
MR {\bf{84j:}}{\rm{46033}}

[S]
Santal\'o, L. A.,
An affine invariant for convex bodies of $n$-dimensional space (Spanish), 
{\it{Portugal. Math.}} {\bf{8}} (1949), 155-161, MR {\bf{12,}}{\rm{526}}



\endRefs

\end